\documentclass{article}
\usepackage{amsmath}
\usepackage{url}

\begin{document}

\title{Construction of the Digamma Function \\by Derivative Definition}
\author{Michael Morales\\
\small\scshape Universidad del Valle de Guatemala\\
\small e-mail: mikael.mm@gmail.com}
\date{}
\maketitle

\begin{abstract}
   The Digamma and Polygamma functions are important tools in mathematical physics, not only for its many properties but also for the applications in statistical mechanics and stellar evolution. In many textbooks is found its develop almost by the same procedure. In this paper expressions for the Digamma and Polygamma functions, in terms of hypergeometric functions, are found through the derivative definition. 
\end{abstract}

\section{Introduction.}
  
  So far, for the definition of the Digamma Function, $\psi(z)$, is usual to find the expression starting from the Euler Infinite Limit definition of $\Gamma(z)$:

\begin{equation}
\Gamma(z)=\lim_{n\rightarrow \infty} \frac{n!}{(z+1)(z+2)...(z+n)}n^z
\end{equation}

 Then taking the natural logarithm and then differentiate, thus obtaining, as appear in \cite{Arfken}, \cite{Ross} and \cite{Chaudry}:

\begin{equation}
\psi(z+1)=\frac{d \ln\Gamma(z+1)}{dz}=-\gamma+\sum^\infty_{n=1} \frac{z}{n(n+z)}
\end{equation}

  where $\gamma$ is the Euler-Mascheroni constant.\\

  However, there's almost nothing about how to obtain $\psi(z)$ from the usual derivative definition

\begin{equation}  \frac{d\Gamma (z)}{dz}=\lim_{h\rightarrow 0} \frac{\Gamma (z+h)-\Gamma(z)}{h} \end{equation}

  In this paper we elaborate an expression of the Digamma function from the derivative definition. It will be shown that is only necessary application of elemental series expansion \cite{Stewart} and some Gamma and Beta functions identities, \cite{Arfken} and \cite{Weisstein},to obtain it. \\

\section{Construction.}

  As mentioned earlier, $ \frac{d ln \Gamma(z)}{dz} $, the Digamma function, will be develop from the derivative definition.  For that, we will need some useful results.\\

  We use the next identities, which can be found in \cite{Arfken} and \cite{Weisstein}:

\begin{equation} 
B(z,h)=\frac{\Gamma(z)\Gamma(h)}{\Gamma(z+h)} 
\end{equation}
\begin{equation}
B(z,h)=\frac{z+h}{h}B(z,h+1)
\end{equation}

  In addition, we use the following Beta function's expansions:

\[ B(z,h)=2\int^\frac{\Pi}{2}_0 \cos^{2z-1}(\theta)\sin^{2h-1}(\theta)d\theta \]
\[ = 2\int^\frac{\Pi}{2}_0\cos^{2z-2}(\theta)\sin^{2h}(\theta)\frac{\cos(\theta)}{\sin(\theta)}d\theta\]
\[ = 2\int^1_0(1-u^2)^{z-1}u^{2h-1}du\]
\[ \approx 2\int^1_0\left[ 1-(z-1)u^2+\frac{(z-1)(z-2)u^4}{2!}-\frac{(z-1)(z-2)(z-3)u^6}{3!}+\cdots\right]u^{2h-1}du\]

  that becomes

\begin{equation} B(z,h) \approx \frac{1}{h}-\frac{(z-1)}{h+1}+\frac{(z-1)(z-2)}{2!(h+2)}-\frac{(z-1)(z-2)(z-3)}{3!(h+3)}+\cdots
\end{equation}

 With the same procedure, we obtain

\begin{equation}
 B(z,h+1)\approx \frac{1}{z}-\frac{h}{z+1}+\frac{h(h-1)}{2!(z+2)}-\frac{h(h-1)(h-2)}{3!(z+3)}+\cdots
\end{equation}

  Starting now from equation (3)

\begin{equation}  \frac{d\Gamma (z)}{dz}=\lim_{h\rightarrow 0} \frac{\Gamma (z+h)-\Gamma(z)}{h} \end{equation}

\[ =\lim_{h\rightarrow 0} \frac{\frac{\Gamma(z)\Gamma(h)}{B(z,h)}-\Gamma (z)}{h} \]
   
\[ =\lim_{h\rightarrow 0} \Gamma(z) \left[\frac{\Gamma(h)-B(z,h)}{h B(z,h)}\right]\]

\begin{equation}
 \frac{1}{\Gamma(z)} \frac{d\Gamma(z)}{dz} = \psi(z) = \lim_{h\rightarrow 0}\frac{\Gamma(h)-B(z,h)}{(z+h)B(z,h+1)} 
\end{equation}

  where $ \psi(z)$  is the Digamma Function and the identities (4) and (5) were used.\\

  We now use Weierstrass's Gamma Function definition for $\Gamma(h)$ and expand

 \[ \lim_{h\rightarrow 0}\Gamma(h) =\lim_{h\rightarrow 0}\frac{1}{h} e^{-\gamma h} \displaystyle\prod^\infty_{n=1} \left( 1+\frac{h}{n}\right) ^{-1}e^{\frac{h}{n}}\]

\[ \approx \lim_{h\rightarrow 0}\frac{1}{h}\left( 1-\gamma h +O(h)^2 \right) \lim_{h\rightarrow 0} \displaystyle\prod^\infty_{n=1} \left( 1+\frac{h}{n}\right) ^{-1}e^{\frac{h}{n}} \]

\begin{equation} \lim_{h\rightarrow 0}\Gamma(h) \approx  \lim_{h\rightarrow 0}\frac{1}{h}-\gamma
\end{equation}

  Substituting the expressions (6), (7) and (10) in equation (9) and taking the limit of h, the desired result is finally obtained

\begin{equation}
\psi(z)=-\gamma+\sum^\infty_{n=1} (-1)^{n+1}\frac{\displaystyle\prod_{i=1}^n(z-i)}{n\, n!}
\end{equation}

\section{Yet another expression for Digamma and Polygamma functions}

   It can be shown that

\[
\displaystyle\prod_{i=1}^n(z-i)=-(-1)^{n+1} \frac{\Gamma(n+1-z)}{\Gamma(1-z)}
\]

\[
=-(-1)^{n+1}(1-z)_n
\]

  where $(1-z)_n$ is the Pochhammer symbol, \cite[eq. 89.1.1]{Hansen}. Using this result in equation (11) leads to

\begin{equation}
\psi(z)=-\gamma-\sum^\infty_{n=1} \frac{(1-z)_n}{{n\, n!}}
\end{equation}

  and in terms of hypergeometric functions

\begin{equation}
\psi(z)=-\gamma-(1-z)\, _3F_2(2-z,1,1;2,2;1)
\end{equation}
    
  Finally we can define the Polygamma Functions as:

  \begin{equation}
\psi^{(l)}(z)=\frac{d\,^l }{dz^l}\, \left( (z-1)_3F_2(2-z,1,1;2,2;1)\right)
\end{equation}

  with $l=1,2,3,\ldots$, for the derivatives of the Digamma function.\\

  As a final comment, for evaluations in the form\\

  \[
\psi_m(z)=-\gamma-\sum^m_{n=1} \frac{(1-z)_n}{{n\, n!}}
\]

  is necessary that $m\geq \Re(z)$ to converge.
    
\section{Conclusion.}    

  It was demonstrated that is possible to find an expression of the digamma function developed from the derivative definition. These expressions, equations (12) to (14), are alternative definitions for the digamma and polygamma functions.

\end{document}